\documentclass{birkmult}
\usepackage{graphicx}
\usepackage{amssymb}
 \newtheorem{thm}{Theorem}[section]
 
 \newtheorem{lem}[thm]{Lemma}
 \newtheorem{prop}[thm]{Proposition}
 \theoremstyle{definition}
 \newtheorem{defn}[thm]{Definition}
 \theoremstyle{remark}
 \newtheorem{rem}[thm]{Remark}
 
 \numberwithin{equation}{section}

\newcommand{\Comp}{\mathbb C}

\renewcommand{\d}{\mathrm{d}}
\newcommand{\Mul}{\mathrm{Mul}}

\begin{document}
%
%
%
%
%
%
%
%

\title{\Large Algebraic Multiplicity and the Poincar\'{e}
Problem}%

\author{Jinzhi Lei}

\address{Zhou Pei-Yuan Center for Applied Mathematics\\
Tsinghua University\\
Beijing, 100084\\
P.R.China}

\email{jzlei@mail.tsinghua.edu.cn}

\thanks{This work was supported
by the National Natural Science Foundation of China(10301006)}

\author{Lijun Yang}

\address{Department of Mathematical Science\\
Tsinghua University\\
Beijing, 100084\\
P.R.China}

\email{lyang@math.tsinghua.edu.cn}

\subjclass{34A05; 34M99}

\keywords{Polynomial system, invariant algebraic curve,
Poincar\'{e} problem}

\date{}

\begin{abstract}
In this paper we derive an upper bound for the degree of the
strict invariant algebraic curve of a polynomial system in the
complex project plane under generic condition. The results are
obtained through the algebraic multiplicities of the system at the
singular points. A method for computing the algebraic multiplicity
using Newton polygon is also presented.
\end{abstract}

\maketitle


\section{Introduction}
\label{sec:1} In this paper, we will present an approach to
establish the upper bound of the degree of the strict invariant
algebraic curve of a polynomial system in the complex projective
plane $\mathbb{P}_\Comp^2$. A polynomial system in
$\mathbb{P}_\Comp^2$ is defined by the vector field
\begin{equation}
\label{eq:27} \dot{z} = P(z,w),\ \ \dot{w} = Q(z,w),
\end{equation}
where  $P$ and $Q$ are relatively prime polynomials with complex
coefficients.

\begin{defn}
\label{def:2} A polynomial $f(z,w)$ is said to be a \index{Darboux
polynomial}Darboux Polynomial of (\ref{eq:27}) if there exists a
polynomial $R_f(z,w)$ such that
\begin{equation}
\label{eq:28} P(z,w)\dfrac{\partial f}{\partial z}+Q(z,w)
\dfrac{\partial f}{\partial w} = R_f(z,w)f(z,w).
\end{equation}
We call the zero-set $C(f)=\{(z,w)\in\hat{\Comp}^2|\,f(z,w)=0\}$
an \index{invariant algebraic curve}invariant algebraic curve, and
$R_f$ the \index{cofactor}cofactor of $f$. In particular, if
$C(f)$ contains no constant irreducible component (i.e., the line
$z = z_0$ or $w = w_0$), then $f$ is a strict Darboux polynomial,
and $C(f)$ is a strict invariant algebraic curve.
\end{defn}

The study of invariant algebraic curves of a polynomial system
goes back to \index{Darboux, Jean-Gaston} Darboux and
\index{Poincar\'{e}, Henri}Poincar\'{e}( see Schlomiuk
\cite{Sch}). In general, the Darboux polynomial of the sytem
\eqref{eq:27} can be found by solving the equation \eqref{eq:28}
for $f$ and $R_f$. Equation \eqref{eq:28} is easy to solve if the
degree of $f$ is known in advance (for example, see
\index{Pereira}Pereira \cite[Propostion 1]{Per:2}). However, it is
still an open problem, for a given system, to establish the upper
bound for the degree of the invariant algebraic curve effectively.
This is named as the \index{Poincar\'{e} problem}Poincar\'{e}
problem. It is known that such an upper bound do exists for a
given polynomial system, see Schlomiuk \cite[Corollary 3.1]{Sch}.
However, the uniform upper bound that depends merely on the degree
of the system does not exists, for non-trivial example, see
Ollagnier \cite{Ol:01}. As consequence, the practical arithmetic
to find the bound from the coefficients is significant for the
general solution for finding the invariant algebraic curve of a
polynomial system. For more remarks and results on the
Poincar\'{e} problem, see Carnicer \cite{Car:94}, Campillo and
Carnicer \cite{Car:97}, Schlomiuk \cite{Sch}, Walcher \cite{Wal}.
The first result to address the Poincar\'{e} problem was presented
by \index{Carnicer, Manuel M.} Carnicer\cite{Car:94} as following.

\begin{thm}[\index{Carnicer's theorem} Carnicer's theorem\cite{Car:94}]
Let $\mathcal{F}$ be a foliation of $\mathbb{P}_{\mathbb{C}}^2$
and let $C$ be an algebraic curve in $\mathbb{P}_{\mathbb{C}}^2$.
Suppose that $C$ is invariant by $\mathcal{F}$ and there are no
dicritical singularities of $\mathcal{F}$ in $C$. Then
$$\partial^o C\leq \partial^o \mathcal{F} + 2$$
\end{thm}

In the proof of Carnicer's theorem, the relationship between the
sum of the multiplicities of a foliation along the branches of a
curve, the degree of the curve, the degree of the foliation and
the Euler characteristic of the curve are systematic used. This
idea is also used in the present paper. However, it was not
provided in \cite{Car:94} the effective method to determine
whether a singular point is dicritical or not. The same inequality
had been showed by Cerveau and Lins Neto\cite{Cer:91} for system
of which all the singularities of the invariant algebraic curve
are nodal. A more straightforward result was presented by
\index{Walcher, Sebastian} Walcher using elementary
method\cite{Wal}. Walcher's result stated:

\begin{thm}{\cite[Theorem 3.4]{Wal}}
\label{th:1} Assume that a vector field $X$ of degree $M$ on
$\mathbb{P}_\Comp^2$ admits an irreducible invariant algebraic
curve, and if all the stationary points of $X$ at infinity are
nondegenerate and non-dicritical, then the degree of the curve
cannot exceed $M+1$.
\end{thm}

In Walcher's proof, the \index{Poincar\'{e}-Dulac normal
form}Poincar\'{e}-Dulac normal forms of the nondegenerate
stationary points of a vector field were discussed. In particular,
when the stationary point is non-dicritical, the precise
information of the number of irreducible semi-invariants of the
vector field $X$ was obtained, from which the upper bound of the
degree of an invariant algebraic curve is derived. It was also
pointed out in \cite{Wal} that if there are dicritical ones among
the nondegenerate stationary points, then the vector field can
admit infinitely many (pairwise relatively prime) semi-invariants.
Moreover, the condition of non-dicritical can be verified through
the investigation of the linear approximation of the vector field
at the stationary points. Thus, Walcher's result provided a
practical approach for the Poincar\'{e} problem.

In this paper, we will present an alternative approach for the
Poincar\'{e} problem by considering the algebraic multiplicities
(see Definition \ref{def:1}) of the singular points of the system,
and obtain an approximate inequality for the upper bound for the
degrees under some generic conditions. The main results of this
paper are:

\begin{thm}
\label{th:2} Consider the differential equation
\begin{equation}
\label{eq:17} \dfrac{\d w}{\d z} = \dfrac{P(z,w)}{z\,Q(z,w)},
\end{equation}
of degree $M = \max\{\deg P(z,w), \deg z\,Q(z,w)\}$. If
\eqref{eq:17} admits an irreducible strict Darboux polynomial
$f(z,w)$. Let $a_1, \cdots, a_k\in \mathbb{C}$ be all roots of
$P(0,w) = 0$, and $a_0 = \infty$, and $\Mul(0,a_i)$ be the
algebraic multiplicity of $(0,a_i)$, then
\begin{equation}
\label{eq:20} \deg_w f(z,w) \leq \sum_{i = 0}^k\Mul(0,a_i).
\end{equation}
In particular, if the singularities $(0,a_i)$ are not algebraic
critical, then
\begin{equation}
\label{eq:P1} \deg_w f(z,w)\leq M\,(k + 1).
\end{equation}
\end{thm}

\begin{thm}
\label{th:main} Consider the polynomial system \eqref{eq:27} of
degree $M = \max\{\deg P(z,w),\\ \deg Q(z,w)\}$, if \eqref{eq:27}
has an invariant straight line $L$, and the singular points at $L$
are not algebraic critical, and if \eqref{eq:27} admits an
irreducible strict Darboux polynomial $f(z,w)$, then
$$\deg f(z,w) \leq M(M+1).$$
\end{thm}

Note that, in Theorem \ref{th:main}, we don't need the
singularities to be non-degenerate, and we will see in next
section that not algebraic critical is weaker than non-dicritical.
In Theorem \ref{th:main}, we require that (\ref{eq:27}) has an
invariant straight line. In fact, it is generic that the line at
infinity is invariant. Hence, the condition in Theorem
\ref{th:main} is generic.

The rest of this paper is arranged as following. In Section
\ref{sec:5}, we will introduce the concept and computing method of
algebraic multiplicity. And next, the main theorems are proved. In
Section \ref{sec:2}, as application, the 2D Lotka-Volterra system
is studied.

\section{Algebraic Multiplicity and the Poincar\'{e} Problem}
\label{sec:5} Let $f(z,w)$ be a Darboux polynomial of
(\ref{eq:27}). In general, the upper bound of the degree of
$f(z,w)$ can not be determined merely from the equation
(\ref{eq:28}). The assumption that $f(z,w)$ is irreducible must be
taken into account. If $f(z,w)$ is irreducible, non lost the
generality, perform the transformation $(z,w)\mapsto (z + c\,w,w)\
(c\in \mathbb{R})$ if necessary, we may assume that $\deg_w f(z,w)
= \deg f(z,w)$. Let $m = \deg_w f(z,w)$, then there are $m$
algebraic functions $w_i(z)$ satisfying $f(z,w_i(z)) = 0\ (i =
1,2,\cdots,m)$. If these $m$ algebraic functions pass through some
common singular points, then $m$ can be bounded by the possible
number of the algebraic solutions that pass through these singular
points. To this end, we will define the algebraic multiplicity as
the number of local algebraic solutions as following.
\begin{defn}
\label{def:1} Consider a differential equation
\begin{equation}\label{eq:16}
\dfrac{\d  w}{\d  z} = F(z,w),
\end{equation}
and $(z_0, w_0)$ $\in \Comp^2$. A formal series
\begin{equation}
\label{eq:30}
w(z) = w_0 + \sum_{i\geq
0}\alpha_i\,(z-z_0)^{\mu_i},\\
\end{equation}
is said to be a \index{local algebraic solution}local algebraic
solution of \eqref{eq:16} at $(z_0,w_0)$ if $w(z)$ is a formal
series solution of \eqref{eq:16} with $\alpha_i\not= 0$,
$\mu_i\in\mathbb{Q}^+$, and $\mu_i < \mu_{i+1} \ (\forall i)$. The
\index{algebraic multiplicity}algebraic multiplicity of
\eqref{eq:16} at $(z_0,w_0)$, denoted by $\Mul(z_0,w_0; F)$ or
simply by $\Mul(z_0,w_0)$ while the context is clear, is defined
as the number of distinct local non-constant algebraic solutions
of \eqref{eq:16} at $(z_0,w_0)$. If $\Mul(z_0,w_0) = \infty$, then
$(z_0,w_0)$ is said \index{algebraic critical}algebraic critical.
\end{defn}

It is evident that algebraic critical implies
\index{dicritical}dicritical (i.e., there are infinitely many
invariant curves passing through the same point).

When $w_0 =\infty$, let $\bar{w} = 1/w$, then $\bar{w}(z)$
satisfies
\begin{equation} \label{eq:32} \dfrac{\d  \bar{w}}{\d  z} =
-\bar{w}^2\,F(z,1/\bar{w}): = \bar{F}(z, \bar{w}),
\end{equation}
and the algebraic multiplicity $\Mul(z_0, \infty; F)$ is simply
defined as $\Mul(z_0,0;\bar{F})$.

Let $a, b, c \in \Comp$ with $a, c\not= 0$, and let $W = a\,(w -
w_0) + b\, (z - z_0)$, $Z=c\,(z - z_0)$, then $W(Z)$ satisfies an
equation of form
\begin{equation}
\label{eq:33} \dfrac{\d  W}{\d  Z} = \tilde{F}(Z,W).
\end{equation}
It is easy to show that a local algebraic solution of
\eqref{eq:16} at $(z_0, w_0)$ corresponds to a local algebraic
solution of \eqref{eq:33} at $(0,0)$. Hence we have
\begin{equation}
\label{eq:6} \Mul(z_0, w_0; F) = \left\{\begin{array}{ll}
\Mul(0, 0; \tilde{F}),& \mathrm{if}\ \tilde{F}(Z,0) \not\equiv 0 \\
\Mul(0, 0; \tilde{F}) + 1,& \mathrm{if}\ \tilde{F}(Z,0) \equiv 0
\end{array}\right.
\end{equation}

It is evident that, if $(z_0, w_0)$ is a regular point and
$F(z,w_0)\not\equiv 0$, then $\Mul(z_0, w_0) = 1$. To estimate the
algebraic multiplicity at singular point $(z_0, w_0)$, we can
substitute \eqref{eq:30} into \eqref{eq:16} to find out all
possible formal series solutions. A method for finding the formal
series solution of a polynomial system at a singular point is
given in Lei and Guan \cite{Lei} using \index{Newton
polygon}Newton polygon (Bruno\cite{Bruno}, Chebotarev\cite{Ceb}).
The result and proof are restated below.

\begin{lem}
\label{le:1}Consider the polynomial system \begin{equation}
\label{eq:22} \dfrac{\d w}{\d z} = \dfrac{P(z,w)}{Q(z,w)}
\end{equation}
 where
$$P(z,w) = \sum_{i \geq 0} P_i(z)\, w^i,\ \ \ Q(z,w) =
\sum_{i \geq 0} Q_i(z)\, w^i,
$$
and
$$ P_i(z) = p_{i,0}\, z^{k_i} +
p_{i,1}\,z^{k_i + 1} + \cdots, \ \ Q_i(z) = q_{i,0}\,z^{l_i} +
q_{i,1}\, z^{l_i + 1} + \cdots,\ \ (i\geq 0)
$$
If $(0,0)$ is a singular point of \eqref{eq:22}, and there exists
$j$, satisfying
\begin{enumerate}
\item[(1).] $k_j = l_{j-1} - 1$; \item[(2).] For any $i\not=j$,
$$\min\{k_i,l_{i-1}-1\} > k_j + (j-i)\,(p_{j,0}/q_{j-1,0})
$$
 \item[(3).] $p_{j,0}/q_{j-1,0}\in
\mathbb{Q}^+$,
\end{enumerate}
then $(0,0)$ is algebraic critical for the system \eqref{eq:22}.
\end{lem}
\begin{proof}
Let $\lambda = p_{j,0}/q_{j-1,0}$, and $u(z) =
w(z)\,z^{-\lambda}$, then $u(z)$ satisfies
\begin{eqnarray*}
\dfrac{\d u}{\d z} &=& \dfrac{\sum_{i\geq 0} (p_{i,0}\,z^{k_i +
i\,\lambda} - q_{i-1,0}\,\lambda\,z^{l_{i-1} + i\,\lambda - 1}+
h.o.t.)\,u^i}{\sum_{i\geq 0}(q_{i,0}\,z^{l_i + (i+1)\,\lambda} +
h.o.t.)u^i}\\
&=&\dfrac{z^{l_{j-1} + j\,\lambda-1}\,\sum_{i\geq 0}
(p_{i,0}\,z^{k_i-k_j + (i-j)\,\lambda} -
q_{i-1,0}\,\lambda\,z^{l_{i-1}-l_{j-1} + (i-j)\,\lambda}+
h.o.t.)\,u^i}{z^{l_{j-1} + j\,\lambda}\,\sum_{i\geq
0}(q_{i,0}\,z^{l_i -l_{j-1} + (i-j)\,\lambda} + h.o.t.)u^i}
 \end{eqnarray*} Taking the conditions of $j$ into account , we can
rewrite above equation as
$$
\frac{\d u}{\d z} =\dfrac{z^{s}\,\hat{P}(z,u)}{z\,\hat{Q}(z,u)}
$$
where $\hat{P}(0,u), \hat{Q}(0,u)\not\equiv 0$, and $s  =
\min_{i\geq 0}\{k_i-k_j + (i-j)\,\lambda, l_{i-1}-l_{j-1} +
(i-j)\,\lambda\}\in \mathbb{Q}^+$. Let $z = \bar{z}^{q_{j-1,0}}$,
then
\begin{equation}
\label{eq:29} \dfrac{\d u}{\d \bar{z}}=
\dfrac{q_{j-1,0}\,\bar{z}^{s\,q_{j-1,0} -
1}\,\hat{P}(\bar{z}^{q_{j-1,0}},u)}{\hat{Q}(\bar{z}^{q_{j-1,0}},u)}
\end{equation}
It's easy to have $s\,q_{j-1,0}\in \mathbb{N}$ and
$\hat{P}(\bar{z}^{q_{j-1,0}},u), \hat{Q}(\bar{z}^{q_{j-1,0}},u)$
are polynomials of $\bar{z}$ and $u$. Thus, for any $\alpha$ such
that $\hat{Q}(0,\alpha)\not=0$, \eqref{eq:29} has a unique
solution $u(\bar{z};\alpha)$ which is analytic at $\bar{z} = 0$
and satisfies $u(0;\alpha) = \alpha$. Thus,
$$w(z;\alpha) =
z^{\lambda}\,u(z^{1/q_{j-1,0}};\alpha) = z^{\lambda}(\alpha +
\sum_{i\geq 1}\frac{1}{i!}u^{(i)}_z(0;\alpha)z^{i/q_{j-1,0}})$$ is
a solution of \eqref{eq:22}, i.e., $w(z;\alpha)$ is a local
algebraic solution of \eqref{eq:22} for any $\alpha$ such that
$\hat{Q}(0;\alpha)\not=0$. Hence, $(0,0)$ is algebraic critical
for \eqref{eq:22}.
\end{proof}
\begin{rem}
\label{re:3}\begin{enumerate} \item The Lemma \ref{le:1} is also
valid for equations of which $P$ and $Q$ are Puiseux series of $z$
and $w$ (with slight change in the proof):
$$P(z,w) = \sum_{i,j\geq0}p_{i,j}\,z^{i/\mu}\,w^{j/\nu},\ \ \ Q(z,w) = \sum_{i,j\geq 0}q_{i,j}\,z^{i/\mu}\,w^{j/\nu}\ \ \ (\mu,\nu\in \mathbb{N})$$
\item From the proof of Lemma \ref{le:1}, if the index $j$
satisfied conditions (1), (2), but
$p_{j,0}/q_{j-1,0}\in\mathbb{R}^+\backslash\mathbb{Q}^+$, let
$\lambda = p_{j,0}/q_{j-1,0}$, then \eqref{eq:22} admits infinity
solutions of form $w(z;\alpha) = z^{\lambda}\,u(z^{1/s};\alpha)$,
where $u(\bar{z};\alpha)$ is the solution of
$$\dfrac{\d u}{\d \bar{z}} = \dfrac{\hat{P}(\bar{z}^{1/s},u)}{s\,\hat{Q}(\bar{z}^{1/s},u)}$$
such that $u(0;\alpha) = \alpha$. Thus, \eqref{eq:22} is
dicritical at $(0,0)$, but not necessary algebraic critical.
\end{enumerate}
\end{rem}

\begin{lem}
\label{th:3} Let $(0,0)$ be a singular point of \eqref{eq:22},
then either $(0,0)$ is algebraic critical, or
\begin{equation}
\label{eq:24} \Mul(0,0)\leq \max\{\deg_w P(z,w), \deg_w Q(z,w) +
1\}.
\end{equation}
\end{lem}
\begin{proof}
Let $N = \deg_w P(z,w), M = \deg_w Q(z,w)$, and
$$P(z,w) = \sum_{i = 0}^N P_i(z)\, w^i,\ \ \ Q(z,w) =
\sum_{i = 0}^M Q_i(z)\, w^i,
$$
where
$$ P_i(z) = p_{i,0}\, z^{k_i} +
p_{i,1}\,z^{k_i + 1} + \cdots, \ \ Q_i(z) = q_{i,0}\,z^{l_i} +
q_{i,1}\, z^{l_i + 1} + \cdots
$$
Substitute
\begin{equation}
\label{eq:w}
w(z) = \alpha_0\, z^{\lambda_0} + h.o.t.\ \ (\alpha_0\not=0, \lambda_0\in\mathbb{Q}^+)
\end{equation}
into \eqref{eq:22}, then \begin{eqnarray*} 0&=&\sum_{i =
0}^MQ_i(z)(\alpha_0\,z^{\lambda_0} +
h.o.t.)^i\,(\alpha_0\lambda_0\,z^{\lambda_0-1} + h.o.t.)
-\sum_{i=0}^NP_i(z)\,(\alpha_0\,z^{\lambda_0} + h.o.t.)^i \\
&=&\sum_{i=0}^M q_{i,0}\,{\lambda_0}\,\alpha_0^{i+1}\,z^{l_i +
(i+1)\,\lambda_0 - 1} - \sum_{i=0}^N p_{i,0}\,\alpha_0^i\,z^{k_i +
i\,\lambda_0} + h.o.t.
\end{eqnarray*}
Thus, at least two of the exponents:$$l_i + (i+1)\,\lambda_0 - 1,
\ \ k_j + j\,\lambda_0,\ \ (0\leq i\leq M,\ \ 0\leq j\leq N)$$are
equal to each other and not larger than any other exponents, and
$\alpha_0\not=0$ that vanishes the coefficient of the lowest
degree. If this is the case, $(\lambda_0, \alpha_0)$ is said to be
acceptable to \eqref{eq:22}. Assume that $(0,0)$ is not algebraic
critical (i.e., Lemma \ref{le:1} is not satisfied), then the
values $\lambda_0$ and $\alpha_0$ can be obtained using Newton
polygon\cite{Bruno, Ceb} as following. Let $\Gamma$ be the Newton
open polygon of all points(see Figure \ref{fig:1})
\begin{equation}
\label{eq:21} (i+1, l_i - 1),\ \ \ (j, k_j),\ \ (0\leq i\leq M,\ \
0\leq j\leq N)
\end{equation}
Let $\Gamma_{i_1}^{i_2}$ be an edge of $\Gamma$, with $i_1 < i_2$
to be the horizontal coordinates of the extreme vertices. Let
$-\lambda_0$ to be the slope of $\Gamma_{i_1}^{i_2}$, then
$\alpha_0$ should satisfy a polynomial of degree $i_2-i_1$. In
particular, $(\lambda_0, \alpha_0)$ is said to be $d$-folded if
$\alpha_0$ is a $d$-folded root of above polynomial. Thus, for the
edge $\Gamma_{i_1}^{i_2}$, there are at most $i_2 - i_1$ pairs of
$(\lambda_0, \alpha_0)$  that are acceptable to \eqref{eq:22}.
Thus, there are totally at most $\max\{M+1,N\}$ pairs of
$(\lambda_0, \alpha_0)$ that are acceptable to \eqref{eq:22}.

\begin{center} \unitlength=0.5cm
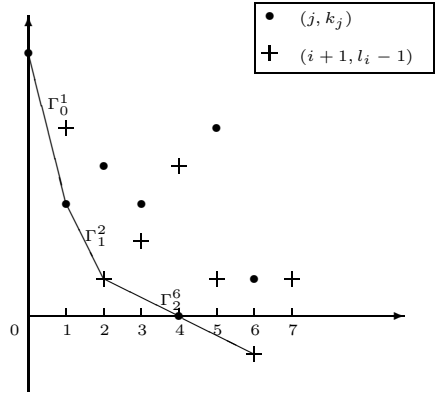
\begin{figure}[hbt]
\begin{picture}(12,12)
\put(0,2){\vector(1,0){10.0}} \put(0,0){\vector(0,1){10.0}}
\put(-0.50,1.5){\tiny{0}} \put(0,9){\circle*{0.2}}
\put(1,5){\circle*{0.2}} \put(1,6.8){\line(0,1){0.4}}
\put(0.8,7){\line(1,0){0.4}} \put(2,2.8){\line(0,1){0.4}}
\put(1.8,3){\line(1,0){0.4}} \put(3,5){\circle*{0.2}}
\put(3,3.8){\line(0,1){0.4}} \put(2.8,4){\line(1,0){0.4}}
\put(4,2){\circle*{0.2}} \put(4,5.8){\line(0,1){0.4}}
\put(3.8,6){\line(1,0){0.4}} \put(5,7){\circle*{0.2}}
\put(5,2.8){\line(0,1){0.4}} \put(4.8,3){\line(1,0){0.4}}
\put(6,0.8){\line(0,1){0.4}} \put(5.8,1){\line(1,0){0.4}}
\put(6,3){\circle*{0.2}} \put(7,2.8){\line(0,1){0.4}}
\put(6.8,3){\line(1,0){0.4}} \put(2.0,6){\circle*{0.2}}
\put(0,9){\line(1,-4){1.0}} 
\put(1,5){\line(1,-2){1.0}} \put(2,3){\line(2,-1){4.0}}
\put(1,2){\line(0,1){0.2}} \put(2,2){\line(0,1){0.2}}
\put(3,2){\line(0,1){0.2}} \put(4,2){\line(0,1){0.2}}
\put(5,2){\line(0,1){0.2}} \put(6,2){\line(0,1){0.2}}
\put(7,2){\line(0,1){0.2}} \put(0.9,1.5){\tiny{1}}
\put(1.9,1.5){\tiny{2}} \put(2.9,1.5){\tiny{3}}
\put(3.9,1.5){\tiny{4}} \put(4.9,1.5){\tiny{5}}
\put(5.9,1.5){\tiny{6}} \put(6.9,1.5){\tiny{7}}
\put(6,10){\framebox[2.4cm][l]{\put(0.2,0){\circle*{0.2}}\ \ \
\put(1.0,-0.2){\tiny{$(j,k_j)$}} \put(0.2,-1.2){\line(0,1){0.4}}
\put(0,-1){\line(1,0){0.4}}\ \ \ \put(1.0,-1.2){\tiny{$(i +
1,l_i-1)$}}}} \put(0.5,7.5){\tiny{$\Gamma_0^1$}}
\put(1.5,4.0){\tiny{$\Gamma_1^2$}}
\put(3.5,2.3){\tiny{$\Gamma_2^6$}}
\end{picture} \caption{Newton Polygon}\label{fig:1}
\end{figure}
\end{center}

For each $(\lambda_0, \alpha_0)$ in the first step, let $w(z) =
\alpha_0\, z^{\lambda_0} + w_1(z)$, then $w_1(z)$ satisfies the
equation
\begin{equation}
\label{eq:w1}
Q(z,\alpha_0\,z^{\lambda_0}+w_1) (\alpha_0\, {\lambda_0}\,
z^{\lambda_0-1}+w'_1)-P(z,\alpha_0\, z^{\lambda_0}+w_1)=0.
\end{equation}
Repeat the foregoing argument, if $(0,0)$ is not algebraic
critical point of \eqref{eq:22}, then there are finite solutions
of \eqref{eq:w1} of form
\begin{equation}
\label{eq:w2}
w_1(z) = \alpha_1\,z^{\lambda_1} + h.o.t. ,\ \ (\lambda_1 \in \mathbb{Q}^+,\ \
\lambda_1 > \lambda_0).
\end{equation}
To complete the proof, it's sufficient to show that if
$(\lambda_0, \alpha_0)$ is $d-$folded, then there are at most $d$
pairs of $(\lambda_1, \alpha_1)$ with $\lambda_1 > \lambda_0$
which are acceptable to \eqref{eq:w1}.

Let \begin{eqnarray*} Q_1(z,w_1) &=& Q(z,\alpha_0\,
z^{\lambda_0} + w_1),\\
P_1(z,w_1) &=& P(z,\alpha_0\, z^{\lambda_0} + w_1) -
\alpha_0\,\lambda_0\,z^{\lambda_0 -
1}\,Q(z,\alpha_0\,z^{\lambda_0} + w_1)
\end{eqnarray*}
then $w_1(z)$ satisfies
\begin{equation}
\label{eq:w11} Q_1(z,w_1)\,w_1' - P_1(z,w_1) = 0
\end{equation}
 Write
$$Q_1(z,w_1) = \sum_{i \geq 0}Q_{1,i}(z)\,w_1^i,\ \ \
P_1(z,w_1) = \sum_{i \geq 0}P_{1,i}(z)\,w_1^i$$ and let $l_{1,i}$
and $k_{1,i}$ be the lowest degrees of $Q_{1,i}(z)$ and
$P_{1,i}(z)$ respectively, and $r_{1,i} =
\min\{k_{1,i},l_{1,i-1}-1\}$. We will prove that if $(\lambda_0,
\alpha_0)$ is $d$-folded, then for any $i>d$,
\begin{equation}
\label{eq:31} r_{1,d} \leq r_{1,i} + (i - d)\,\lambda_0
\end{equation}
When \eqref{eq:31} is satisfied, then there are at most $d$-pairs
of $(\lambda_1,\alpha_1)$ which are acceptable to \eqref{eq:w11}
and $\lambda_1 > \lambda_0$. In fact, let $(\lambda_1,\alpha_1)$
to be acceptable to \eqref{eq:w11}, then there exist $j_1 < j_2$,
such that
$$\lambda_1 = \dfrac{r_{1,j_1} - r_{1,j_2}}{j_2 - j_1} > \lambda_0$$
and
$$r_{1,d} \geq r_{1,j_1} + (j_1 - d)\,\lambda_1,\ \ r_{1,d} \geq r_{1,j_2} + (j_2 - d)\,\lambda_1 $$
If $j_1 > d$ (or $j_2 > d$), then
$$r_{1,d} >  r_{1,j_1} + (j_1 -
d)\,\lambda_0\ \ \ (\mathrm{or}\ \ r_{1,d} >  r_{1,j_2} + (j_2 -
d)\,\lambda_0)$$ which is contradict to \eqref{eq:31}. Hence,
$j_1<j_2 \leq d$, and there are at most $d$-pairs of
$(\lambda_1,\alpha_1)$ (taking account that $(0,0)$ is not
algebraic critical).

To prove \eqref{eq:31}, let
\begin{eqnarray*}
&&Q(z,\alpha\, z^{\lambda_0}) = \sum_{i \geq 0}\xi_i(\alpha)\,
z^{s_i}\ \ \ \ \ \ \ \ \ (s_0 < s_1 < \cdots) \\
&&P(z,\alpha\, z^{\lambda_0}) = \sum_{i \geq 0}\eta_i(\alpha)\,
z^{\tau_i}\ \ \ \  \ \ \ \ \ (\tau_0 < \tau_1 <\cdots)\\
\end{eqnarray*}
then
\begin{eqnarray*}
Q_{1,i}(z)
&=&\dfrac{1}{i!}\,z^{-i\,\lambda_0}\,\sum_{j \geq
0}\xi_j^{(i)}(\alpha_0)\, z^{s_j}\\
P_{1,i}(z)
&=&\dfrac{1}{i!}\,z^{-i\,\lambda_0}\,\left(\sum_{j\geq
0}\eta_j^{(i)}(\alpha_0)\,z^{\tau_j} -
\alpha_0\,\lambda_0\,z^{\lambda_0 - 1}\,\sum_{j \geq 0}
\xi_j^{(i)}(\alpha_0)\,z^{s_j} \right)
\end{eqnarray*}
and hence
\begin{equation}
\label{eq:46}r_{1,i} \geq \min\{\tau_0, s_0 + \lambda_0 - 1\}  -
i\, \lambda_0.
\end{equation}
Thus, it is sufficient to show that
\begin{equation} \label{eq:23} \min\{k_{1,d}, l_{1,d-1} - 1\} =
\min\{\tau_0, s_0 + \lambda_0 -1\} - d\,\lambda_0.
\end{equation}
To this end, write
\begin{eqnarray*}
Q_{1,d-1}(z) &=& \frac{1}{d!}
\xi_0^{(d-1)}(\alpha_0)\,z^{s_0 + \lambda_0 - d_0\,\lambda_0} + h.o.t.\\
P_{1,d}(z) &=&
\frac{1}{d!}\left(\eta_0^{(d)}(\alpha_0)\,z^{\tau_0} -
\alpha_0\,\lambda_0\,\xi_0^{(d)}(\alpha_0)\,z^{s_0 + \lambda_0 -
1}\right)\cdot z^{-d\,\lambda_0} + h.o.t.
\end{eqnarray*}
and let
$$P(z,\alpha\,z^{\lambda_0})  - \alpha\,\lambda_0\,z^{\lambda_0 - 1}\,Q(z,\alpha\,z^{\lambda_0}) =  \varphi(\alpha)\,z^{v_0} + h.o.t.$$
Because $(\lambda_0,\alpha_0)$ is acceptable to \eqref{eq:22} and
$d$-folded, we have
\begin{equation}
\label{eq:f} \varphi(\alpha_0) = \cdots =
\varphi^{(d-1)}(\alpha_0) = 0,\ \varphi^{(d)}(\alpha_0)\not = 0.
\end{equation}
Therefore, we have the following:
\begin{enumerate}
\item[(a).] If $\tau_0 < s_0 + \lambda_0-1$, then $\varphi(\alpha)
= \eta_0(\alpha)$ and $\eta_0^{(d)}(\alpha_0) \not = 0$.
\item[(b).] If $s_0 + \lambda_0 -1 < \tau_0$, then
$\varphi(\alpha) = -\lambda_0\,\alpha\, \xi_0(\alpha)$, and hence
$\xi_0^{(d)}(\alpha_0) \not = 0$. \item[(c).] If $s_0 + \lambda_0
-1 = \tau_0$, then $\varphi_0(\alpha) = \eta_0(\alpha) - \alpha
\lambda_0 \xi_0(\alpha)$, and hence
$$\varphi_0^{(d)}(\alpha_0)=- \lambda_0 \xi_0^{(d-1)}(\alpha_0) + (\eta_0^{(d)}(\alpha_0)  - \alpha_0 \lambda_0
\xi_0^{(d)}(\alpha_0)) \not = 0.$$ Thus, we have
$\xi_0^{(d-1)}(\alpha_0)\not=0$ or $\eta_0^{(d)}(\alpha_0) -
\alpha_0 \lambda_0\xi_0^{(d)}(\alpha_0)\not=0$.
\end{enumerate}

It is not difficult to verify that \eqref{eq:23} is held in any
one of the above cases, and thus the Lemma is concluded.
\end{proof}

From the proof of Lemma \ref{th:3}, the local algebraic solutions
of \eqref{eq:22} at $(0,0)$ can be obtained by repeating the
Newton polygon. Moreover, following the procedure, we will either
stop by the case that $(0,0)$ is algebraic critical (Lemma
\ref{le:1}), or encounter the local algebraic solution of form
$$w(z) = \sum_{i = 0}^{k-1}\alpha_i\,z^{\lambda_i} + u(z)$$
where $(\lambda_{k-1}, \alpha_{k-1})$ is 1-folded, and $u(z)$
satisfies an equation
\begin{equation}
\label{eq:38} \dfrac{\d u}{\d z} =
\dfrac{\hat{P}(z,u)}{\hat{Q}(z,u)}
\end{equation}
 where $\hat{P},
\hat{Q}$ are Puiseux series. Whenever this is the case, we have
the following.
\begin{lem}
\label{le:2} In the equation \eqref{eq:38} that derived from
\eqref{eq:22} through above procedure, let
$$\hat{P}(z,u) = \hat{p}_{0,0}z^{k_0} + \hat{p}_{1,0}z^{k_1}\,u + h.o.t.,
\ \ \hat{Q}(z,u) = \hat{q}_{0,0}z^{l_0} + h.o.t.$$ If
$(\lambda_{k-1}, \alpha_{k-1})$ is 1-folded, and one of the
following is satisfied:
\begin{enumerate}
\item[(1).]$k_1\not=l_0-1$; or \item[(2).] $k_1=l_0-1$, and
$\hat{p}_{1,0}/\hat{q}_{0,0}\not\in (\lambda_{k-1},\infty)\cap
\mathbb{Q}^+$,
\end{enumerate}
then $(0,0)$ is not algebraic critical of \eqref{eq:22}.
\end{lem}
\begin{proof}
Let $u(z)$ be a local algebraic solution of \eqref{eq:38},
expressed as
\begin{equation} \label{eq:43} u(z) = \sum_{i\geq
k}\alpha_i\,z^{\lambda_i}
\end{equation}
where $\lambda_i > \lambda_{i-1},\ (\forall i\geq k)$. We will
show that $(\lambda_i, \alpha_i)$ are determined by \eqref{eq:38}
uniquely.

From the proof of Lemma \ref{th:3},  we have
$$k_0 - \min\{k_1,l_0-1\} > \lambda_{k-1}$$
Hence, substitute \eqref{eq:43} into \eqref{eq:38}, and taking
account that $(\lambda_{k-1}, \alpha_{k-1})$ is 1-folded, and
either $k_1\not=l_0-1$ or $k_1=l_0-1$, $p_{1,0}/q_{0,0}\not\in
(\lambda_{k-1},\infty)\cap \mathbb{Q}^+$, we have $\lambda_k = k_0
- \min\{k_1,l_0-1\}$, and $\alpha_k$ is determined uniquely by
$p_{0,0}, q_{0,0}, p_{1,0}, k_1,l_0$. Therefore, $(\lambda_k,
\alpha_k)$ is also 1-folded. Let $u(z) = \alpha_k\,z^{\lambda_k} +
v(z)$, then $v(z)$ satisfies
\begin{equation}
\label{eq:42}\dfrac{\d v}{\d z} = \dfrac{\hat{p}'_{0,0}\,z^{k_0'}
+ \hat{p}_{1,0}\,z^{k_1}\,v + h.o.t.}{\hat{q}_{0,0}\,z^{l_0} +
h.o.t.}
\end{equation}
 where $k_0' > k_0$. In
particular, conditions in the Lemma are also valid for
\eqref{eq:42}. Thus, we can repeat the procedure, and hence there
is unique solution $u(z)$ of form \eqref{eq:43}, and $(0,0)$ is
not algebraic critical for \eqref{eq:22}.
\end{proof}
\begin{rem}
\label{re:1} In the Lemma \ref{le:2}, we might also find the
solution of form \eqref{eq:43} when $k_1=l_0-1$ and
$\hat{p}_{1,0}/\hat{q}_{0,0}\in (\lambda_{k-1},\infty)\cap
\mathbb{Q}^+$. However, when this is the case, we can identify two
cases:
\begin{enumerate}
\item[(1).] If $\hat{p}_{1,0}/\hat{q}_{0,0} \in (\lambda_i,
\lambda_{i+1})\cap \mathbb{Q}$ for some $i\geq k-1$, then the
condition in Lemma \ref{le:1} is satisfied at the $i$'th step, and
$(0,0)$ is algebraic critical. \item[(2).] If
$\hat{p}_{1,0}/\hat{q}_{0,0}  = \lambda_i$ for some $i$, then
$(0,0)$ is not algebraic critical.
\end{enumerate}
In any case, we can stop the procedure in finite steps. Thus, it's
effective to find the algebraic multiplicities of \eqref{eq:22}
using the Newton polygon.
\end{rem}

\textbf{Example}\ \ Consider the equation
\begin{equation}
\label{eq:45} (z + w^2)\,w' - (z^2 + \mu w) = 0
\end{equation}

\begin{center} \unitlength=0.5cm
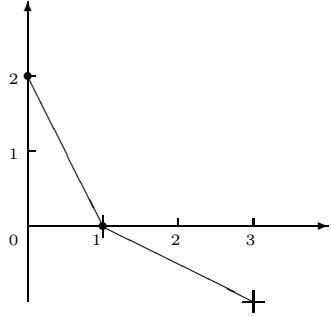
\begin{figure}[hbt]
\begin{picture}(8,8)
\put(0,2){\vector(1,0){8.0}} \put(0,0){\vector(0,1){8.0}}
\put(0,6){\circle*{0.2}} \put(2,2){\circle*{0.2}}
\put(1.7,2.0){\line(1,0){0.6}} \put(2.0,1.7){\line(0,1){0.6}}
\put(5.7,0){\line(1,0){0.6}} \put(6.0,-0.3){\line(0,1){0.6}}
\put(0,6){\line(1,-2){2.0}} \put(2,2){\line(2,-1){4.0}}
\put(2,2){\line(0,1){0.2}} \put(4,2){\line(0,1){0.2}}
\put(6,2){\line(0,1){0.2}} \put(-0.5,1.5){\tiny{0}}
\put(1.7,1.5){\tiny{1}} \put(3.8,1.5){\tiny{2}}
\put(5.8,1.5){\tiny{3}} \put(0,4){\line(1,0){0.2}}
\put(0,6){\line(1,0){0.2}} \put(-0.5,3.8){\tiny{1}}
\put(-0.5,5.8){\tiny{2}}
\end{picture} \caption{Newton Polygon of (\ref{eq:45})}\label{fig:2}
\end{figure}
\end{center}

The Newton polygon of \eqref{eq:45} is shown at Figure
\ref{fig:2}. From the Newton polygon, if $\mu\in (1/2, 2)\cap
\mathbb{Q}$, then $(0,0)$ is algebraic critical, with local
algebraic solutions
$$w(z) = \alpha_0 z^{\mu} + h.o.t.\ \ \ (\alpha_0\not=0)$$
Mean while, if $\mu\not\in (1/2, 2)\cap \mathbb{Q}$, the possible
local algebraic solutions are
\begin{eqnarray*}
w(z) &=& \frac{1}{2-\mu}\,z^2 + h.o.t.\ \ (\mathrm{if}\ \
\mu\not=2)\\
w(z) &=&\pm\sqrt{2\mu-1}\,z^{1/2} + h.o.t.\ \ (\mathrm{if}\
\mu\not=1/2)
\end{eqnarray*}
When $\mu\not=2$, let
$$w(z) = \frac{1}{2-\mu}\,z^2 + w_{1,1}(z) $$
then $w_{1,1}(z)$ satisfies
$$w_{1,1}' = \dfrac{2\,z^5 - (2-\mu)^3\,\mu w_{1,1} + h.o.t.}{-(2-\mu)^3\,z + h.o.t.}$$
Thus, we conclude the following. If $\mu\in (2,5)\cap \mathbb{Q}$,
then $(0,0)$ is algebraic critical, with local algebraic solutions
$$w(z) = \frac{1}{2-\mu}\,z^2 + \alpha_1\,z^\mu + h.o.t,\ \ \ (\alpha_1\not=0).$$
If $\mu\not=2,5$, we have the local algebraic solution
$$w(z) = \frac{1}{2-\mu}\,z^2 - \frac{2}{(5-\mu)\,(2-\mu)^3}\,z^5 + h.o.t.$$
When $\mu\not\in [1/2,2)\cap\mathbb{Q}$, let
$$w(z) = \sqrt{2\mu-1}\,z^{1/2} + w_{1,2}(z)$$
then $w_{1,2}(z)$ satisfies
$$w_{1,2}' = \dfrac{2\,z^{5/2} + (2-2\mu)\,z^{1/2}\,w_{1,2} + h.o.t.}{4\,\mu\,z^{3/2} + h.o.t.}$$
Thus, if $\mu\not=1/5$, we have the local algebraic solution
$$w(z) = \sqrt{2\mu-1}\,z^{1/2} + \frac{1}{5\,\mu - 1}\,z^2 +h.o.t.$$

Thus, repeat the above procedure, we can determine, for given
$\mu$, the algebraic multiplicity $\Mul(0,0)$ of \eqref{eq:45}. In
particular, if $\mu\not\in (1/2,\infty)\cap \mathbb{Q}$, then
$\Mul(0,0)\leq 3$. \begin{flushright} $\square$\end{flushright}

In the rest of this section, we will prove the main results.

\begin{proof}[Proof of Theorem \ref{th:2}]
Let $W$ be the set of all non-constant local algebraic solutions
of \eqref{eq:17} at $(0,a_i)$ for some $0\leq i\leq k$. Then
$$|W| = \sum_{i=0}^k\Mul(0,a_i)$$
Let $f(z,w)$ be an irreducible strict Darboux polynomial of
\eqref{eq:17}, and $m = \deg_wf(z,w)$, then there are $m$
algebraic functions $w_i(z)$ that defined by $f(z,w) = 0$. It is
sufficient to show that any algebraic function $w_i(z)\in W$. To
this end, we only need to show that
\begin{equation}
\label{eq:44}\lim_{z\to 0}w_i(z)= \{a_0,a_1,\cdots,a_k\}
\end{equation}
Consider the equation
$$z\,Q(z,w)\,\dfrac{\partial f}{\partial z} +
P(z,w)\,\dfrac{\partial f}{\partial w} = R_f(z,w)\,f(z,w)$$ Let $z
= 0$, then $f(0,w)$ satisfies
$$P(0,w)\,f_w'(0,w) = R_f(0,w)\,f(0,w).$$
Thus $f(0,w)$ is an non-constant multiply of $\prod_{i = 1}^k (w -
w_i)^{l_i},\ \ (l_i\geq 0)$. From which \eqref{eq:44} is easy to
conclude.

It is easy to have $\Mul(0,\infty) \leq M$. Hence, if the
singularities $(0,a_i)$ are not algebraic critical, then, from
Lemma \ref{th:3},
$$\deg_wf(z,w)\leq M\,(k+1)$$
\end{proof}

\begin{proof}[Proof of Theorem \ref{th:main}]
If \eqref{eq:27} has an invariant straight line $L$,  perform
suitable transformation, we may assume that $L$ is given by
$$a\,z + b\,w + z = 0,\  \ \ \ (a\not=0)$$
and $\deg f(z,w) = \deg_w f(\frac{z - b\,w - c}{a}, w)$. It is
easy to see that the degree of the system is not increase under
linear transformation. Let
$$\bar{w} = w, \bar{z} = a\,z + b\,w + c,$$
then $\bar{w}(\bar{z})$ satisfies the equation of form
\begin{equation}
\label{eq:18} \dfrac{d \bar{w}}{d \bar{z}} =
\dfrac{\bar{P}(\bar{z}, \bar{w})}{\bar{z}\,\bar{Q}(\bar{z},
\bar{w})},
\end{equation}
where $\bar{P}(\bar{z}, \bar{w}), \bar{Q}(\bar{z}, \bar{w})$ are
polynomials. Moreover, $\bar{f}(\bar{z}, \bar{w}) =
f(\frac{\bar{z} - b\,\bar{w} - c}{a}, \bar{w})$ is an irreducible
Darboux polynomial of \eqref{eq:18}, and $\deg f(z,w) =
\deg_{\bar{w}} \bar{f}(\bar{z}, \bar{w})$. Let $(a_i, b_i), (1\leq
i \leq M)$ be singular points of \eqref{eq:27} at $L$, then $(0,
b_i)$ are singular points of \eqref{eq:18} at $\bar{z} = 0$, and
not algebraic critical. Hence, apply Theorem \ref{th:2} to
\eqref{eq:18}, we have
$$\deg f(z,w) = \deg_{\bar{w}} \bar{f}(\bar{z}, \bar{w}) \leq M\,(M + 1).$$
\end{proof}

\section{Application to 2D Lotka-Volterra system}
\label{sec:2} In this section, we will apply Theorem \ref{th:2} to
2D \index{Lotka-Volterra system}Lotka-Volterra system:
\begin{equation}
\label{lv:2-1} \dot{z} = z\,(z + c\,w - 1),\ \ \ \dot{w} =
w\,(b\,z + w - a).
\end{equation}
Invariant algebraic curves of Lotka-Volterra system had been
studied by many authors. Recent results on this topic, refer to
Ollaginer \cite{Ol:01-2}, Cair\'{o} \textit{et.al.}\cite{JL:03}
and the references. In Ollaginer \cite{Ol:01-2}, the complete list
of parameters of which the system has strict invariant algebraic
curve is presented. We will reobtain one part of the results
through the algebraic multiplicity.

Note that (\ref{lv:2-1}) is invariant under following
transformations:
\begin{eqnarray}
\label{eq:11} (z,w,a,b,c)&\rightarrow& (\frac{w}{a}, \frac{z}{a},
\frac{1}{a}, c, b),\ \mathrm{if}\ a\not=0;\\
\label{eq:12} (z,w,a,b,c)&\rightarrow& (\frac{1}{z},
(1-c)\frac{w}{z}, 1-b, 1-a, \frac{c}{c-1}),\ \mathrm{if}\ c\not=1.
\end{eqnarray}
Results in this section are also valid under above
transformations.

Since  $z = 0$ and $w = 0 $ are invariant straight lines of
\eqref{lv:2-1}, Theorem \ref{th:2} is applicable.
\begin{prop}
\label{main:1} If the 2D L-V system
\begin{equation}
\label{lve} \dfrac{\d  w}{\d  z} =
\dfrac{w\,(b\,z+w-a)}{z\,(z+c\,w-1)}
\end{equation}
has a strict Darboux polynomial $f$, then
$$
\begin{array}{rcll}
\deg_w f(z,w)&\leq& \Mul(0,\infty) + \Mul(0,a) +
\Mul(0,0),&\ \  \mathrm{if}\ \  a\not=0 \\
\deg_w f(z,w)&\leq& \Mul(0,\infty) + \Mul(0,0).&\ \ \mathrm{if}\ \   a = 0 \\
\end{array}
$$
\end{prop}

In particular, we have.
\begin{prop}
\label{main:2} If in \eqref{lve},
\begin{equation}
\label{eq:1} a\not\in \mathbb{Q}^+, c\not\in \mathbb{Q}^-, c -
\dfrac{1}{a}\not\in\mathbb{Q}^+ \backslash\{1\}
\end{equation}
then \eqref{lve} has strict invariant algebraic curve if and only
if
$$a (1-c) + (1-b) = 0,$$
and the invariant algebraic curve is given by
$$a (z-1) + w = 0.$$
\end{prop}
\begin{proof}
When \eqref{eq:1} is satisfied, the singularities $(0,0), (0,a),
(0,\infty)$ are not algebraic critical, and
$$\Mul(0,0) = 0, \Mul(0,a)\leq 1, \Mul(0,\infty) = 0$$

Hence, if $f(z,w)$ is a strict irreducible Darboux polynomial,
then $\deg_w f = 1$. From which the result is easy to conclude.
\end{proof}

Proposition \ref{main:2} shows that the algebraic multiplicities
may give an exact bound for the degree of the Darboux polynomial
in particular cases. However, if there are algebraic critical
points among the singularities, \eqref{eq:20} does not provide the
finite value. In this case, as we had seen from Lemma \ref{le:1},
there are infinite local algebraic solutions. On the other hand,
this does not automatically imply that all these local algebraic
solutions are algebraic solutions. And hence, we come to the
following concrete problem: If a singular point of a system is
algebraic critical, how many local algebraic solutions are exactly
the algebraic function? It requires additional work to discuss
this problem, and one may hope that the solution of this problem
should lead to the final resolution of the Poincar\'{e} problem.



\section*{Acknowledgement}
The authors would like to thank professor Zhiming Zheng and
Dongming Wang for their work at the organization of the seminar of
DESC2004. The authors are also grateful to the referees for their
helpful suggestion.


\end{document}